\documentclass[twocolumn]{revtex4}
\usepackage{times,amsmath,amsfonts,amsthm}

\makeatletter       % Macros for arabic section numbering

\renewcommand{\p@subsection}{}

\renewcommand{\p@subsubsection}{}

\makeatother

\newtheorem{theorem}{Theorem}

\newtheorem{corollary}[theorem]{Corollary}
\newtheorem{remark}{Remark}
\newtheorem*{definition}{Exchangeability Condition}

\def\({\left(}
\def\){\right)}
\def\[{\left[}
\def\]{\right]}
\def\a{\alpha}  \def\d{\delta} 
   
   \def\p{\phi} 
  \def\t{\theta} 
 
\def\Cal{\cal}    \def\pp{{\Cal P}}  \def\xx{{\Cal X}} 
\def\Array{\begin{eqnarray*}}
\def\EndArray{\end{eqnarray*}}
\def\Enumerate{\begin{enumerate}}
\def\EndEnumerate{\end{enumerate}}
\def\Eq{\begin{equation}}
\def\EndEq{\end{equation}}
\def\EqArray{\begin{eqnarray}}
\def\EndEqArray{\end{eqnarray}}
\def\mref#1{(\ref{#1})}
\def\qt#1{\qquad\text{#1}}
\def\as{\stackrel{\mathrm{a.s}}{\rightarrow}}

\def\iid{\stackrel{\mathrm{i.i.d.}}{\sim}}

\def\p{{\bf p}}
\def\Pr{\mathbb{P}}

\begin{document}

\title{Random Probability Measures via P\'olya Sequences:\\
Revisiting the Blackwell-MacQueen Urn Scheme}
\author{Hemant Ishwaran}
\email{ishwaran@bio.ri.ccf.org}
\author{Mahmoud Zarepour}
\email{zarepour@mathstat.uottawa.ca}
\affiliation{Cleveland Clinic Foundation and University of Ottawa}

\begin{abstract}
Sufficient conditions are developed for a class of generalized P\'olya
urn schemes ensuring exchangeability.  The extended class includes
the Blackwell-MacQueen P\'olya urn and the urn schemes for the
two-parameter Poisson-Dirichlet process and finite dimensional
Dirichlet priors among others.
\end{abstract}
\maketitle

\section{Introduction}
By making use of a remarkably simple generalized P\'olya urn scheme,
Blackwell and MacQueen (\cite{BM}) described an elegant alternate way to
approach the Ferguson Dirichlet process~(\cite{FG}).  Let
$X_1,X_2,\ldots$ be a sequence of random elements on a complete
separable metric space $\xx$ defined by
\Eq
\Pr\{X_1\in\cdot\}=\mu(\cdot)/\mu(\xx)
\label{blackwell.one}
\EndEq
and 
\Eq
\Pr\{X_{i+1}\in\cdot|X_1,\ldots,X_i\}=\mu_i(\cdot)/\mu_i(\xx),
\qt{}i\ge 1,
\label{blackwell.two}
\EndEq
where $\mu_i(\cdot)=\mu(\cdot)+\sum_{j=1}^{i}\d_{X_i}(\cdot)$ and
$\mu$ is a finite non-null measure on $\xx$.  Blackwell and MacQueen
called such a sequence a {\it P\'olya sequence with parameter $\mu$}.
They showed that if $X_1,X_2,\ldots$ was such a sequence, then:
\Enumerate
\item[(a)] $\mu_i(\cdot)/\mu_i(\xx)$ converges almost surely to a
discrete random probability measure $\mu^*$.
\item[(b)] $\mu^*$ is the Ferguson Dirichlet process with parameter $\mu$.
\item[(c)] Given $\mu^*$, $X_1,X_2,\ldots$ are independent with 
           distribution $\mu^*$.
\EndEnumerate

Result (c) shows that the Blackwell-MacQueen P\'olya sequence is
exchangeable, while (b) shows that the sequence is an infinite sample
from the Dirichlet process.  Thus (a) and (b) combined show that the
P\'olya urn defined by~\mref{blackwell.one} and~\mref{blackwell.two}
is a way to draw values from the Dirichlet process.  Moreover, (a)
shows that the Dirichlet process is the limit for the urn distribution
$\mu_i(\cdot)/\mu_i(\xx)$, thus providing an alternate way to
characterize the Dirichlet process.  These facts are far more
difficult to prove than the contrapositive result which starts with a
sample from the Dirichlet process and shows that such a sample must be
exchangeable and can be constructed from a P\'olya urn.  The latter
result follows from elementary properties of the Dirichlet process
which we now describe.  Let $X_1,X_2,\ldots$ be a sequence derived
from a Dirichlet process with parameter $\mu$, i.e:
\Array
(X_i|P) &\iid& P,\qt{}i=1,2,\ldots\\
P       &\sim&\mu^*.
\EndArray
It was shown in~\cite{FG} that $\mu^*(\cdot|X_1,\ldots,X_i)$, the
posterior for $\mu^*$ based on the first $i$ observations
$X_1,\ldots,X_i$, is also a Dirichlet process, but with an updated
parameter $\mu_i$ (see~\cite{LW}, Section 3.2, for a proof using
a Laplace functional argument).  An immediate consequence of this is
that
\Array
&&\hskip-20pt\Pr\{X_{i+1}\in\cdot|X_1,\ldots,X_i\}\\
&=&\int\Pr\{X_{i+1}\in\cdot|X_1,\ldots,X_i,P\}\,\mu^*(dP|X_1,\ldots,X_i) \\
&=&\int P(\cdot)\,\mu^*(dP|X_1,\ldots,X_i) \\
&=&\mu_i(\cdot)/\mu_i(\xx),
\EndArray
and, thus, $X_1,X_2,\ldots$ can be defined by the
P\'olya urn described by~\mref{blackwell.one}
and~\mref{blackwell.two} (the fact that the sequence is exchangeable
follows by definition).

Put another way, elementary properties for the Dirichlet process shows
us that the {\it prediction rule}, that is, the conditional
distribution for $X_{i+1}$ given $X_1,\ldots,X_i$, corresponds to an
exchangeable generalized P\'olya urn distribution
$\mu_i(\cdot)/\mu_i(\xx)$.  This type of direct result is somewhat
unique as it is generally hard to derive simple explicit prediction
rules for a general random discrete probability measure.  Instead,
another way to approach the problem is in the direction studied by
Blackwell and MacQueen.  Thus, it is natural to wonder what types of
P\'olya urn schemes other than~\mref{blackwell.one}
and~\mref{blackwell.two} lead to: (i) an exchangeable sequence
$X_1,X_2,\ldots$ and (ii) an urn distribution with a limiting random
discrete probability measure?  Sufficient and necessary conditions for
(i) have been given in~\cite{HP}~(cf Theorem 2) in terms
of what is called the exchangeable partition probability function
(EPPF), a symmetric non-negative function which characterizes the
distribution of an exchangeable partition on the positive integers
$\{1,2,\ldots\}$.  In this paper, however, we take a more direct
approach to answering (i) (and consequently (ii)), by introducing an
Exchangeability Condition (Section 3) which puts constraints on the
manner in which the urn scheme selects new values or chooses a
previously sampled value.  While these conditions will be shown only
to be sufficient to ensure (i), they have the advantage that they are
simpler to understand then conditions stated in terms of the EPPF.
The proof should also be readily accessible to non-experts to this
area.  Our Exchangeability Condition is shown to hold for several
important generalized P\'olya urns, including those for the
two-parameter Poisson-Dirichlet process~(\cite{PY}) as well as the
class of finite dimensional Dirichlet priors~(\cite{IZa}).  Corollary
1 of Section 3, our main result, summarizes our results.

\section{Notation and background}
We begin by introducing some notation necessary to explain our
generalized P\'olya urn schemes.  Let $\p_i=\{C_{j,i}:
j=1,\ldots,n(\p_i)\}$ denote a partition of $\{1,\ldots,i\}$ where
$C_{j,i}$ is the $j$th set of the partition.  Write $e_{j,i}$ for the
cardinality of $C_{j,i}$.  Thus, $\p_i$ is a partition made of $n(\p_i)$
sets and $\sum_{j=1}^{n(\p_i)} e_{j,i}=i$.  Let $X_1^*,X_2^*,\cdots$
denote the sequence of unique values {\it in the order of their
appearance} from $X_1,X_2,\ldots$ and let $\p_i$ be a partition of
$\{1,\cdots,i\}$ recording the clustering of the first $i$
observations $X_1,\ldots,X_i$.  By this we mean $X_l=X_j^*$ for each
$l\in C_{j,i}$, where $j=1,\ldots,n(\p_i)$.

Let $\nu$ denote a non-null probability measure over $\xx$.  We
will consider sequences $X_1,X_2,\ldots$ defined by
\Eq
\Pr\{X_1\in\cdot\}=\nu(\cdot)
\label{general.one}
\EndEq
and 
\EqArray
&&\hskip-20pt\Pr\{X_{i+1}\in\cdot|X_1,\ldots,X_i\}\label{general.two}\\
&=&\frac{q_{0,i}}{\sum_{j=0}^{n(\p_i)}q_{j,i}}\nu(\cdot)+
\sum_{j=1}^{n(\p_i)}\frac{q_{j,i}}{\sum_{j=0}^{n(\p_i)}q_{j,i}}\d_{X_j^*}(\cdot),
\qt{}i\ge 1,\nonumber
\EndEqArray
where 
\Eq
q_{j,i}:= q_{j,i}(e_{1,i},\ldots,e_{n(\p_i),i}),\;
q_{0,i}:= q_{0,i}(e_{1,i},\ldots,e_{n(\p_i),i})\\
\label{qj.def}
\EndEq
are non-negative real valued symmetric functions depending only upon
$\{e_{1,i},\ldots,e_{n(\p_i),i}\}$.

The form for $q_{0,j}$ and $q_{j,i}$ in~\mref{qj.def} is suggested by
Theorem~1 of~\cite{HP} which states that for $X_1,X_2,\cdots$ to be
exchangeable, the functions $q_{0,i}$ and $q_{j,i}$ must be almost
surely equal to some function of the partition $\p_i$; or
equivalently, they must be some function of the cardinalities
$e_{j,i}$.  For example, observe that the Blackwell-MacQueen P\'olya
sequence (with parameter $\mu$) corresponds to the choices
$q_{0,i}=\mu(\xx)$, $q_{j,i}=e_{j,i}$ and
$\nu(\cdot)=\mu(\cdot)/\mu(\xx)$.

In proving our general result, an important technical condition that
we will need to address concerns the choice for $\nu$.  We say that
$\nu$ is {\it non-atomic} if $\nu\{x\}=0$ for each $x\in\xx$.  One of
the unique features of the Blackwell-MacQueen P\'olya urn scheme is
that it yields an exchangeable sequence regardless of whether $\nu$ is
non-atomic.  For example, if $\xx=\{1,\ldots,r\}$ is a finite sample
space and $\mu$ is a finite discrete measure over $\xx$,
then~\mref{blackwell.one} and~\mref{blackwell.two} implies that
$X_1,\ldots,X_i$ is the result of successive draws from an urn
originally having $\mu(l)$ balls of color $l$, and following each draw
for a ball, the ball is replaced and another one of its same color is
added to the urn.  It follows that
$$
\Pr\{X_1=x_1,\ldots,X_i=x_i\} 
=\frac{1}{\mu(\xx)^{[i]}}\prod_{l=1}^r \mu(l)^{[n(l)]},
\label{finite.dp}
$$
where $n(l)$ denotes the number of $x$'s equal to $l$ and
$a^{[i]}=a(a+1)\cdots(a+i-1)$ (note: $a^{[0]}=1$).  Observe that since
the right-hand side is a symmetric function of $(x_1,\ldots,x_i)$, it
follows automatically that $X_1,\ldots,X_i$ is exchangeable.  This
fact was not lost on Blackwell and MacQueen.  Indeed, the key to their
proof relies on the fact that their P\'olya urn scheme produces an
exchangeable sequence for finite sample spaces.

\subsection{Non-exchangeability over discrete spaces}
However, prediction rules for random discrete measures often break
down when $\xx$ is allowed to be a finite sample space.
A good example is the two-parameter Poisson-Dirichlet
process discussed in~\cite{PY}.  This is the random
discrete probability measure whose prediction rule for a non-atomic $\nu$
corresponds to the choices 
$$
q_{0,i}=\t+\a n(\p_i)
\hskip10pt\text{ and }\hskip10pt
q_{j,i}=e_{j,i}-\a,
$$
where $0\le\a<1\text{ and }\t>-\a$.
See~\cite{PY} and also~\cite{P95}, \cite{P96} for further
details.  Setting $\a=0$ and $\t=\mu(\xx)$ leads to the
Blackwell-MacQueen P\'olya sequence with parameter $\mu=\t\nu$, and as
discussed produces an exchangeable sequence without constraint to
$\nu$.  In general, however, if $\a\ne 0$, exchangeability breaks down
if $\xx$ is allowed to be a finite sample space and $\nu$ is atomic.
This can be easily demonstrated by the following counter-example.  Let
$\xx=\{1,\ldots,r\}$ where $r\ge 2$ and suppose that $\nu(l)=1/r$ for
each $l=1,\ldots,r$.  Then,
\Array
&&\hskip-20pt\Pr\{X_1=1,X_2=2,X_3=1\}\\
&=&\Pr\{X_1=1\}\times\Pr\{X_2=2|X_1\}\\
&&\qquad\times\Pr\{X_3=1, X_3=X_1, X_3\ne X_1|X_1,X_2\}\\
&=&\frac{(\t+\a)\Bigl((\t+2\a)/r+1-\a\Bigr)}{r^2(\t+1)(\t+2)}.
\EndArray
Note that the last expression in the middle equation underlies the
problem with working with an atomic measure.  Here the conditional
event that $X_3=1$ occurs if we choose the previous value $X_1$ or if
we choose the value $X=1$ randomly from $\nu$.  This wouldn't be a
problem with a non-atomic probability measure since the probability of
obtaining a previously observed $X_i$ value is zero under $\nu$.  But
this leads to a breakdown of exchangeability for an atomic measure.
Consider the probability,
\Array
&&\hskip-20pt\Pr\{X_1=1,X_2=1,X_3=2\}\\
&=&\Pr\{X_1=1\}\times\Pr\{X_2=1,X_2=X_1,X_2\ne X_1|X_1\}\\
&& \qquad\times\Pr\{X_3=2|X_1,X_2\}\\
&=&\frac{\Bigl((\t+\a)/r+1-\a\Bigr)(\t+\a)}{r^2(\t+1)(\t+2)}.
\EndArray
Thus, $\Pr\{X_1=1,X_2=2,X_3=1\}\ne \Pr\{X_1=1,X_2=1,X_3=2\}$ unless
$\a=0$.  This shows that only the
Blackwell-MacQueen P\'olya urn is exchangeable in general
for the two-parameter process.

\section{Main results}
Thus, given the technical difficulties in working with atomic
measures, we will hereafter restrict attention to non-atomic measures
$\nu$.  Our results will also rely on the following key conditions for
the functions $q_{0,i}$ and $q_{j,i}$ appearing in~\mref{general.one}
and~\mref{general.two}.

\begin{definition}
For each $i\ge 1$, $q_{j,i}=\psi(e_{j,i})$ and $q_{0,i} =
\psi_0(n(\p_i))$, where $\psi$ and $\psi_0$ are some fixed non-negative
real valued functions.  Furthermore, for each partition $\p_i$ of
$\{1,\ldots,i\}$
\Eq
\sum_{j=0}^{n(\p_i)} q_{j,i} = \xi(i)>0
\label{exch.cond}
\EndEq
where $\xi$ is some fixed positive real valued function.
\end{definition}

These conditions are satisfied by many interesting generalized P\'olya
urn schemes, which we now list.  By satisfying the Exchangeability
Condition, Theorem~1 (stated later) shows that each of these urns
(subject to a non-atomic $\nu$) are exchangeable.
\Enumerate
\item
Independent and identically distributed sampling.  This is of
course the most obvious form of exchangeability and follows with
choices $q_{0,i}=1$ and $q_{j,i}=0$.
\item
$N$ values selected at random.  Let $N>1$ be a positive integer and
let $q_{0,i}=(N-n(\p_i))I\{n(\p_i)<N\}$ and $q_{j,i}=1$.  Observe that
$q_{0,i}$ becomes zero when $n(\p_i)\ge N$, which restricts the
urn sequence from having more than $N$ distinct values.  Note that
condition~\mref{exch.cond} is satisfied because $\sum_{j=0}^{n(\p_i)}
q_{j,i} = N$.
\item
The Blackwell-MacQueen P\'olya sequence with parameter $\mu$.  This
corresponds to $q_{j,i}=e_{j,i}$, $q_{0,i}=\mu(\xx)$ and
$\sum_{j=0}^{n(\p_i)} q_{j,i} = \mu(\xx) + i$.
\item
The two-parameter Poisson-Dirichlet process.  As discussed,
this corresponds to choices $q_{0,i}=\t+\a n(\p_i)$ and
$q_{j,i}=e_{j,i}-\a$, where $0\le\a<1$ and $\t>-\a$.  Thus,
$\sum_{j=0}^{n(\p_i)} q_{j,i} = \t+i$.  This is the prediction
rule for the random discrete probability measure $\pp$ defined by
\EqArray
&&\hskip5pt\pp(\cdot)=V_1\d_{Z_1}(\cdot)\label{stickbreak.def}\\
&&\quad+\sum_{k=2}^\infty\left\{(1-V_1)(1-V_2)\cdots(1-V_{k-1})\,V_k\right\}
\,\d_{Z_k}(\cdot),\nonumber
\EndEqArray
where $\{V_k\}$ are i.i.d Beta$(1-\a,\t+k\a)$ random variables,
independent of $\{Z_k\}$, which are i.i.d with law $\nu$. See~\cite{PY}
for details.  Observe that by setting $\a=0$ we end up
with the Dirichlet process with parameter $\mu=\t\nu$.  In this
case,~\mref{stickbreak.def} corresponds to the stick-breaking
representation for the Dirichlet process.  See~\cite{IJ}
for background on stick-breaking priors.
\item
Finite dimensional Dirichlet priors (Fisher's model).  Let $N>1$ be a
positive integer and let $q_{0,i}=\t(1-n(\p_i)/N)I\{n(\p_i)<N\}$ and
$q_{j,i}=e_{j,i}+\t/N$, where $\t>0$.  Then $\sum_{j=0}^{n(\p_i)}
q_{j,i} = \t+i$ which satisfies~\mref{exch.cond}.  Observe
that the choice for $q_{0,i}$ restricts the process from having more
than $N$ distinct values.  One can show that the values $q_{0,i}$ and
$q_{j,i}$ correspond to the prediction rule for the finite dimensional
Dirichlet prior $\pp_N$ defined by
$$
\pp_N(\cdot)=\sum_{k=1}^N \frac{G_k}{\sum_{k=1}^N G_k}
\,\d_{Z_k}(\cdot),
$$
where $\{G_k\}$ are i.i.d Gamma$(\t/N)$ random variables, independent
of $\{Z_k\}$, which are i.i.d with law $\nu$. See~\cite{P95}, \cite{P96}
and~\cite{IJ} for further details.  Also see~\cite{IZb}
who showed that $\pp_N$ is a weak limit
approximation to the Dirichlet process.
\EndEnumerate

\subsection{Exchangeability}
We now show that our Exchangeability Condition is sufficient to ensure
that the sequence defined by~\mref{general.one} and~\mref{general.two}
is exchangeable.

\begin{theorem}
\label{exch.theorem} 
If $\nu$ is a non-atomic (and non-null) probability measure over $\xx$
and the Exchangeability Condition holds, then $X_1,X_2,\ldots$ is
exchangeable.  
\end{theorem}

\begin{proof}
Let $i>1$ (the case $i=1$ is obvious) and let $dx_1,\ldots,dx_i$
denote a sequence of differentials, some of which can be equal.  Let
$\p_i=\{C_{j,i}: j=1,\ldots,n(\p_i)\}$ be the partition of
$\{1,\ldots,i\}$ which records the clustering of $dx_1,\ldots,dx_i$.
That is, if $dx_1^*,\ldots,dx_{n(\p_i)}^*$ denote the unique values of
$dx_1,\ldots,dx_i$, then $dx_l=dx_j^*$ for each $l\in C_{j,i}$.  As
before, write $e_{j,i}$ for the cardinality of $C_{j,i}$.  For
notational convenience set $\psi(0)=1$.

It follows from the assumption that $\nu$ is non-atomic, and upon
using~\mref{general.one} and~\mref{general.two}, that
\EqArray
&&\hskip-20pt\Pr\{X_1\in dx_1,\ldots,X_i\in dx_i\}\label{marginal.distribution}\\
&=& \prod_{j=1}^i \Pr\{X_j\in dx_j|X_1,\ldots,X_{j-1}\}\nonumber\\
&=&D_i^{-1}\prod_{k=1}^{n(\p_i)-1}\hskip-5pt\psi_0(k)
\prod_{j=1}^{n(\p_i)}\biggr\{\nu(dx_j^*)I\{dx_l=dx_j^*:l\in C_{j,i}\}\nonumber\\
&&\qquad\qquad\times\Bigl(\psi(1)\times\cdots\times\psi(e_{j,i}-1)\Bigr)\biggr\},
\nonumber
\EndEqArray
where the first product (in square brackets) follows from the
assumption that $q_{0,i} = \psi_0(n(\p_i))$ (note: if $n(\p_i)=1$ the
product is assumed to be 1), while the second product uses the
assumption that $q_{j,i}=\psi(e_{j,i})$.  The expression $D_i$
appearing in~\mref{marginal.distribution} is a normalizing constant.
By~\mref{exch.cond}, it can be seen that
$D_i=\xi(1)\times\cdots\times\xi(i-1)$.  Thus, deduce that the
right-hand side of~\mref{marginal.distribution} is a symmetric
function of $(dx_1,\ldots,dx_i)$, and hence that $X_1,\ldots,X_i$ is
exchangeable.  
\end{proof}

\begin{remark}
As a special case, the expression~\mref{marginal.distribution} yields
the well known joint density for the Dirichlet process:
\Array
&&\hskip-20pt\Pr\{X_1\in dx_1,\ldots,X_i\in dx_i\}=
\frac{\mu(\xx)^{n(\p_i)}}{\mu(\xx)^{[i]}}\\
&&\times\prod_{j=1}^{n(\p_i)}\biggl(\nu(dx_j^*)
   I\{dx_l=dx_j^*:l\in C_{j,i}\biggr)
   (e_{j,i}-1)!.
\EndArray
(Substitute $\xi(j)=\mu(\xx)+j$, $\psi(j)=j$ and $\psi_0(j)=\mu(\xx)$
into~\mref{marginal.distribution}).
\end{remark}

\subsection{The Blackwell-MacQueen generalization}
By appealing to Proposition~11 of~\cite{P96}, in combination with
Theorem~1, we obtain the following corollary which is a generalization
of the Blackwell and MacQueen result.

\begin{corollary}
Let $X_1,X_2,\ldots$ be the sequence defined by~\mref{general.one}
and~\mref{general.two} where $\nu$ is a non-atomic (non-null) probability
measure over $\xx$ and $\{q_{0,i},q_{j,i}\}$ satisfy the
Exchangeability Condition.  
\Enumerate
\item[(a)]
Let $F_{i+1}$ denote the conditional distribution for $X_{i+1}$
defined by~\mref{general.two}.  Then, $F_{i+1}\as\pp^*$ in
$\ell_1$-distance, where $\pp^*$ is the random probability measure
defined by
$$
\pp^*(\cdot)=\sum_j p_j\,\d_{X_j^*}(\cdot) +
(1-\sum_j p_j)\nu(\cdot) ,
$$
where $p_j=\lim_{i\rightarrow\infty}e_{j,i}/i$.
\item[(b)] 
$\{X_j^*\}$ are i.i.d $\nu$ and independent of $\{p_j\}$.
\item[(c)]
Given $\pp^*$, $X_1,X_2,\ldots$ are independent with distribution $\pp^*$.
\item[(d)]
If $q_{0,i}/\xi(i)\as 0$, then $\pp^*$ is discrete with probability
one; i.e. $\pp^*(\cdot)=\sum_j p_j\,\d_{X_j^*}(\cdot)$.  
\EndEnumerate
\end{corollary}

\begin{proof} Theorem~1 ensures that $X_1,X_2,\ldots$ is exchangeable.  Thus,
(a), (b) and (c) follows from de Finetti's representation for
exchangeable sequences.  See Theorem~6 of~\cite{P95} and
Proposition~11 of~\cite{P96}.  To prove (d) we use a theorem of~\cite{KP}
which states that if $X_1,X_2,\ldots$ is an
exchangeable sequence from a random measure $\pp^*$, then $\pp^*$ is
discrete with probability one if 
$$
a_i=\Pr\{X_{i+1}\text{ is different than }X_1,\ldots,X_i\}
\as 0 .
$$
See~\cite{SC}, Section 1.6, for a proof.  Thus, (d) is proven
since $a_i=q_{0,i}/\xi(i)$.
\end{proof}

\begin{remark} 
A little bit of work shows that each of our examples listed earlier
(excluding our first example for the i.i.d case) are examples of
generalized P\'olya urn schemes which satisfy condition (d).  Thus,
each produce exchangeable sequences from a random discrete probability
measure.  
\end{remark}

\end{document}